\documentclass[11pt]{amsart}
\usepackage{latexsym,amsfonts,amssymb,amsmath}
\usepackage{pb-diagram}
\usepackage{graphics}

\usepackage[all]{xy}

\newtheorem{theorem}{Theorem}[section]
\newtheorem{lemma}[theorem]{Lemma}
\newtheorem{proposition}[theorem]{Proposition}
\newtheorem{corollary}[theorem]{Corollary}

\newtheorem{algorithm}[theorem]{Algorithm}
\newtheorem{remark}[theorem]{Remark}
\newtheorem{Example}[theorem]{Example}

\newenvironment{example}{\begin{Example}\rm}{\end{Example}}

\def\CC{{\mathbb C}}
\def\NN{{\mathbb N}}
\def\PP{{\mathbb P}}
\def\QQ{{\mathbb Q}}
\def\RR{{\mathbb R}}
\def\ZZ{{\mathbb Z}}

\def\Maple{{\sc Maple\; }}
\def\Singular{{\sc Singular \; }}

\begin{document}

\title{Some results on inhomogeneous discriminants}

\author{Mar\'{\i}a Ang\'elica Cueto and Alicia Dickenstein}
\address{Departamento de Matem\'atica\\
FCEN, Universidad de Buenos Aires \\(1428) Buenos Aires, Argentina.}
\email{alidick@dm.uba.ar, macueto@gmail.com}
\thanks{Both authors are partially supported by UBACYT X042, CONICET
PIP 5617 and ANPCyT PICT03-20569}

\begin{abstract}
We study generalized Horn-Kapranov rational parame\-tri\-za\-tions of inhomogeneous
sparse discriminants from  both a theoretical and an algorithmic perspective. We show
that all these discriminantal parametrizations are birational, we prove some results 
on the corresponding implicit equations, and we propose a combinatorial algorithm to
compute their degree in the uniform case of (co)dimension $3$. 
\end{abstract}

\maketitle

\section{Introduction}

Given a configuration of $n$ lattice points $A$ in $\ZZ^{d-1}$, let
${ F_A = \sum_{a \in A} x_a t^a }$  denote the generic polynomial
in $d-1$ variables $(t_1, \dots, t_{d-1})$ with
exponents in $A$. Under certain general conditions, Gelfand, Kapranov
and Zelevinsky \cite{gkzbook} showed that there exists an
irreducible polynomial with integer coefficients $D_A = D_A(x)\, $
in the vector of coefficients $x = (x_a : a \in A)$ (defined up to
sign),
called the $A$-discriminant, which vanishes for each choice
of coefficients $c$ for which $F_A$ and all its partial
derivatives have a common root in the torus $(\CC^*)^{d-1}$.
The $A$-discriminant
is an affine invariant, in the sense that any configuration
of points affinely isomorphic to $A$ has the same discriminant.

Given $A$, we form the $d \times n$ integer matrix (also called $A$)
whose first row consists of ones, and whose columns are given by
the points $(1,a)$ for all $a \in A$. The kernel of this matrix
expresses the affine dependencies among the given configuration of
points.
Let $B=(b_{ij}) \in \ZZ^{n \times (n-d)}$ be a matrix whose column
vectors are a
basis of the integer kernel of the matrix $A$, i.e. a {\em Gale dual
\/} of $A$.
$B$ is of full rank and its maximal minors have  g.c.d. $g_B=1$.
Since the first row of $A$ is the vector $(1,\dots,1)\in \ZZ^n$,
the row vectors of $B$ must always add up to $0$. We say
that $B$ is {\em regular} when this condition is satisfied.

Set $m = n-d$. The $A$-discriminant $D_A$ is $A$-homogeneous,
i.e it is quasi-homogenous relative to the
weight defined by any vector in the row span of $A$. Therefore,
``taking out these homogeneities''  we get a polynomial
$\Delta_B$ in $m$ variables $y_1, \dots,  y_m$ which is in
fact the implicit equation of (the closure of) the image  $S_B$ of
the rational Horn-Kapranov parametrization  \cite{gkzbook,
log_gauss_map, dfs} defined
by
\begin{equation} \label{eq:psiB}
y_j \, = \, \prod_{i=1}^n (b_{i1} u_1 + \dots +b_{im} u_m)^{b_{ij}},\quad j=1,
\dots, m.
\end{equation}
Extracting the homogeneities means the following.  As we
said, there exists $v \in \ZZ^d$ such that all monomials
$c^\nu$ that occur in $D_A = \sum_{\nu}  d_\nu  \, x^{\nu}$
satisfy $A \cdot \nu \, = \, v$. Therefore, for
any $\nu_0$ such that $A \cdot \nu_0 \, = \, v$,
$$D_A (x) \, = \,  x^{\nu_0} \, \, \sum_{\nu} d_\nu \, x^{\nu - \nu_0},$$
where $d_\nu \in \ZZ \smallsetminus \{0\}$ and $\nu - \nu_0 \in
\ker_\ZZ(A)$.  Write each $\nu - \nu_0$ as a 
linear combination
of the columns $v^{(1)}, \dots, v^{(m)}$ of $B$.
Then, there is a Laurent polynomial $\Delta_B(y)$ in $m$ variables such
that up to a monomial, $\Delta_B (x^{v^{(1)}}, \dots, x^{v^{(m)}})$
equals $D_A(x)$. In particular, $\Delta_B$ has the same number of
monomials and the same coefficients as $D_A$.

\smallskip

Via the well known Cayley trick, the computation of mixed sparse
resultants may be reduced to the computation of sparse discriminants~\cite{gkzbook, rhf}.
The computation of  sparse discriminants $D_A$
(including all sparse resultants),
 is then equivalent
to the implicitization problem for the parametric varieties
given by (\ref{eq:psiB}). These implicitization problems are
thus interesting and very hard, since they involve the whole
of sparse elimination.

\smallskip

It is easy to see that when the matrix $B \in \ZZ^{n \times m}$
defining the parametrization
(\ref{eq:psiB}) is not of maximal rank $m$, then $S_B$ is certainly
not a hypersurface. This is the case in which the corresponding homogenized
discriminant $D_A$ is just the constant polynomial $1$.
However, the image of the parametrization may fail to be a hypersurface
even if $B$ is of maximal rank (c.f. Example~\ref{ex:def}).
The classification of such {\em defective} configurations is a hard combinatorial
problem \cite{cc, dfs, gkzbook}. On the other side, the condition of
non defectiveness can be
checked algorithmically.

When trying to compute discriminants, natural reductions as in
\cite{cc} lead to integer matrices $B$ which are still regular and of
maximal rank, but whose columns do not necessarily generate a
saturated lattice $\ZZ B$, i.e. such that the gcd of its maximal minors
takes any value $g_B \in \ZZ \setminus \{0\}$.
On the other side, sparse discriminants describe the singular locus of $A$-hypergeometric
systems  of PDE's, while the inhomogeneous discriminants $\Delta_B$
describe the singular locus of classical Horn hypergeometric differential
equations \cite{gkz89}. In this setting, it is again  natural and necessary
to consider matrices $B$ with arbitrary gcd $g_B$ \cite{dms05}. Another interest
on studying this case stems in the precise version we give of a result of Kapranov 
\cite{log_gauss_map}, 
which asserts that the rational hypersurfaces 
they define have an interesting geometric characterization.

\smallskip

In this article we focus then on rational parametrizations of the form
(\ref{eq:psiB}) for regular non defective integer matrices $B \in
\ZZ^{n \times m}$ of
maximal rank $m$ and any value of gcd $g_B \not= 0$. To emphasize the fact
that the lattice generated by the columns of the matrix
is not necessarily saturated  (equivalently, the lattice generated by
the rows is
not necessarily $\ZZ^m$), we will call our matrix
$C$ and we will keep the name $B$ in case $g_B=1$. The irreducible equation of
the closure of the image $S_C$, denoted by $\Delta_C \in
\ZZ[y_1,\dots, y_m]$, will
be referred to as a generalized inhomogeneous discriminant.

We prove in Section~\ref{sec:one}
that, up to multiplication by an element $\lambda \in (\CC^*)^{m}$, the
parametrizations \`a la Horn-Kapranov associated to such matrices $C$
correspond precisely to
those rational hypersurfaces $S_C$ for which the logarithmic Gauss map
is birational.
In fact, we follow line by line Kapranov's proof in \cite{log_gauss_map} and we
correct slightly his statement in the sense that the condition $g_C=1$ is not
necessary. In particular, we show that all such parametrizations are birational,
thus recovering and explaining the results in \cite[Section~3.3]{passare}.

We then study in Section~\ref{sec:two} the precise relation between
the generalized
inhomogeneous discriminant $\Delta_C$ which gives the
equation of $S_C$ and the equation $\Delta_B$ of $S_B$, where the
columns of $B$ are a $\ZZ$-basis
of the saturated of the lattice $\ZZ C$ spanned by the columns of $C$. We also give other tips to simplify
the search for the equation $\Delta_C$.

Finally, in Section~\ref{sec:three} we focus on the case $m=3$.  The
case of codimension
$m=2$ has been studied in detail in \cite{DicStu02_elim_2}. We
give a combinatorial algorithm for the
computation of the degree of $S_C$ in case $C$ is uniform, based on
the intersection
formula~(\ref{eq:fulton})  and the algorithm for the
computation of
local multiplicites at the base points given in \cite[Prop.~1.5]{vol_alg}.
This degree could be also obtained via Theorem~2.2 in \cite{dfs}
which should be made explicit to give a dehomogenized version of Theorem~1.3 in that
paper. It would
be interesting to develop intersection formulas which do not require to take
a common denominator in order to work in projective space, to avoid the new
base points coming from this approach. We also
give several examples which
show the difficulty of computing the local multiplicities in the general
case.

\section{Generalized inhomogeneous discriminants and birational Gauss maps}
\label{sec:one}

We start by setting some notations and first consequences of the de\-fi\-ni\-tion 
of the Horn-Kapranov rational parametrizations $\psi_C$ associated to $C$-matrices 
and their closed images $S_C$. 
We then show that these parametrizations are proper (i.e. ${\rm deg} (\psi_C) =1$) 
and that (torus translates of)  varieties of the form $S_C$  
give all rational varieties with birational logarithmic Gauss map.

\subsection{The setting}\label{ss:setting}

Given a  matrix $C \in \ZZ^{n \times m }$, $n \geq m$,
we denote by $C_1, \dots, C_n \in {\ZZ^m}$ the row vectors
of $C$.
Each $C_k$ defines a linear form
\begin{equation}
 {l_{k}}{(u_1,\ldots, u_{m})} := \; \langle C_k \,,\, (u_1,\ldots,
u_{m}) \rangle.
\end{equation}

Throughout this paper, we will always assume that $C$ is regular
(i.e. its columns sum up to zero), and has no zero rows.
We  associate to the matrix $C$ an algebraic variety $S_C$ in $\CC^m$ 
which is the closure of the image of the rational parametrization:
\[
\psi_C: \CC^m \dashrightarrow \CC^m
\qquad (u_1, \ldots, u_{m}) \mapsto (y_1, \ldots, y_m) \, ,
\]
where 
\begin{equation}
y_k  \, = \, \prod_{i=1}^{n}{l_{i}(u_1, \ldots, u_{m})}^{c_{i,k}} 
\qquad \forall  k =1, \ldots, m.
\end{equation}
Call $f_0 \, = \, \prod_{i = 1}^{n}{l_{i}^{-min\{0,\, c_{i,k}\,:\, k =1, \ldots,
  m\}}} $ the least common denominator of all the $y_k$'s, and write
\begin{equation} \label{eq:yf}
y_k \, = \, \frac{f_k}{f_0},  \quad k=1, \ldots m.
\end{equation}
Since the coordinates of the parametrization are given as a product of
linear forms whose exponents sum up to $0$, all $f_i$ have the same degree
$d_C$, where
\begin{equation} \label{eq:d}
d_C \,  =  \, 
 -\sum_{i=1}^n{\min\{0,c_{i,k}\,:\, k =1,
\ldots, m\}}.
\end{equation}
We can think of the mapping
$\psi_C$ as a rational function between projective spaces
\begin{equation} \label{eq:ppsi}
\psi_C: \PP^{m-1} \dashrightarrow \PP^m,
\end{equation}
where $\psi_C \, = \, (f_0 : f_1 : \dots : f_m)$ is defined outside
the {\em base point locus} ${\mathcal Z}\, = \, V(f_0, \dots, f_m)$.
We again denote by $S_C$ the projective variety defined by 
the closure of the image of this map.

The linear forms $l_k$ define a hyperplane arrangement in $\CC^m$ and in
$\PP^{m-1}$. Let $\mathcal F$ be a flat in this arrangement, i.e. a linear space
defined by the vanishing of a subset  $\{l_{i_1}, \dots, l_{i_r}\}$ of the given
linear forms. Denote by ${\mathcal L} (\mathcal F)$  all linear forms $l_j$ vanishing
on $\mathcal F$, i.e. all linear forms $l_j$ which lie in the linear span of
$\{l_{i_1}, \dots, l_{i_r}\}$.  Note that a given $f_k$ vanishes on $\mathcal F$ if and only if 
it contains a linear factor from ${\mathcal L}(\mathcal F)$. A given flat ${\mathcal F}$ of
the arrangement will be called {\em basic} if all of $f_0, \dots, f_m$ vanish on ${\mathcal F}$.
We then have

\begin{lemma} \label{lemma:bpl}
The base point locus $\mathcal Z$ equals the union of all basic flats.
\end{lemma}

\begin{remark} \label{rmk:finitebp}
By taking out any common factor from $f_0, \dots, f_m$, 
we can always assume that ${\rm codim}({\mathcal Z}) \geq 2$.
Then, in the case of codimension $3$ that we will study in Section~\ref{sec:three} we can
assume without loss of generality that the number of base points is finite. 
Note however that their local structure can be very complicated.
When $m >3$, the base point locus variety has in general positive dimension. 
\end{remark}

Recall that $C$ is called defective if ${\rm codim} (S_C) > 1$.
We easily have

\begin{lemma}\label{lemma:def}
Let $C \in \ZZ^{n \times m }$ be a regular matrix of rank $r <m$. 
Then ${\rm codim}(S_C) >  m-r$, so $C$ is defective.
\end{lemma}

\begin{proof} 
Assume $C_1, \dots, C_r$ are linearly independent and write any other
$C_j = \sum_{i=1}^r \lambda_{i,k}^j C_k, \, j=r+1, \dots, n$. It follows
that $y_1, \dots, y_m$ can be written as homogeneous rational functions of 
$l_1, \dots, l_r$ of degree $0$, and so they can be written as rational
functions of the $r-1$ variables $l_1/l_r, \dots, l_{r-1}/l_r$. Therefore,
the codimension of $S_C$ is at least $m -(r-1) = m-r +1$, as claimed.
\end{proof}

\begin{example} \label{ex:def}
As we have remarked in the Introduction, the converse to Lemma~\ref{lemma:def} is not
true and the classification of defective matrices involves subtle combinatorial
questions. As a simple example of a full rank defective matrix, let $n = 2 n'$ be
even and $m \leq n'$ arbitrary. Pick any set of integer vectors $C_1, \dots C_{n'}$
whose $\QQ$-linear span is $\QQ^m$ and take $C_{n'+k} = - C_k$ for all $k =1, \dots, n'$.
Then, it is straightforward to check that $\psi_C$ is a constant map.
\end{example}


We will assume from now on that $C$ has full rank $m$, so that the gcd of its
maximal minors $g_C$ is non zero.

\subsection{Birationality of $\psi_C$}

Given an algebraic group $G$ with Lie algebra $\mathcal{G}$,  let
$l_g:h \mapsto gh$ denote the left translation by an element $g \in G$. 
Let $S   \subset G$ be an irreducible algebraic hypersurface. The left
\emph{Gauss map} of the hypersurface $S$ is the rational map
$\gamma_S: S \rightarrow \PP(\mathcal{G}^*)$, taking a smooth point
$y \in S$ to the hyperplane $d(l_y^{-1})(T_yS)\subset T_e(G)=\mathcal{G}$.
The case which we are interested in is $G={(\CC^*)}^m$. So
$\mathcal{G}=\CC^m$, and left translation is the usual coordinatewise
multiplication map by an element in  the torus. 
If $S$ is a hypersurface defined by a minimal equation
$(\Delta =0)$, then the Gauss map is simply
  \begin{equation}
    \label{eq:Gauss_map}
    \gamma(y)=(y_1 \frac{\partial{\Delta}}{\partial{y_1}}({y}): 
\ldots:  y_m  \frac{\partial{\Delta}}{\partial{y_m}}({y}))\;,
\end{equation}
mapping a regular point ${y} \in S$ to a projective point in $\PP^{m-1}$.
Assume $S^* = S \cap (\CC^*)^m$ is non empty. Geometrically, the map $\gamma_S$
corresponds to looking at the image of $S^*$ via the map $\log(y) = 
(\log(y_1), \dots, \log(y_m))$
and then considering the Gauss map of $\log(S^*)$.

The following theorem is essentially due to Kapranov~\cite{log_gauss_map}. Our
contribution lies in making the statement precise by removing the incorrect hy\-po\-the\-sis
about the gcd of the matrix defining the parametrization. We will show in 
Section~\ref{sec:two}
why a quick first thought may lead to conjecture that this gcd should equal $1$
(cf. Remark~\ref{rmk:cazabobos}).

\begin{theorem}
\label{th:Kapranov}
Let $S \subset \CC^m$ be an algebraic irreducible hypersurface. 
\par
The Gauss map $\gamma_S: S \dashrightarrow
\PP^{m-1}$ is birational if and only if there exist a non defective and regular
integer matrix
$C \in \ZZ^{n \times m}$ of full rank, and a constant $\lambda \in (\CC^*)^m$ 
such that $S= \lambda \cdot S_C$, i.e. $S$ is a torus translate by $\lambda$ of a generalized
inhomogeneous discriminant hypersurface.  
\par
Moreover, in this case, $\lambda \cdot \psi_C$ is birational and 
the logarithmic Gauss map $\gamma_S$ is its inverse.
\end{theorem}

As we said in the Introduction, the proof is exactly Kapranov's original proof, except that one
can check that the hypothesis $g_C =1$ is not used in the ``if'' direction and that  in the
``only if'' direction his last statement about the fact that $g_C$ needs to be equal to $1$ for
$\psi_C$ to be of degree $1$, is not true. For the convenience of the reader, we sketch only the ``if''
direction, i.e. the proof that for any $C$ as in the statement, the logarithmic Gauss map $\gamma_{S_C}$
is indeed a birational inverse to $\psi_C$.

\smallskip

We first give a simple example.

\begin{example}[The discriminant of a generic
univariate cubic polynomial] \label{cubic}
The $A$-discriminant associated to the $2 \times 4$ matrix
$$ A \, = \, \left(
\begin{array}{cccc}
1 & 1 & 1 & 1 \\
0 & 1 & 2 & 3
\end{array}
\right).$$
i.e. the discriminant $D_A(x_0, x_1, x_2,
x_3)$ of the generic polynomial $F_A = x_0 + x_1 t + x_2 t^2 +
x_3 t^3$ equals
\[D_A(x) \, = \, -27 x_3^2 x_0^2+18 x_3 x_0 x_2 x_1+x_2^2 x_1^2-4 x_2^3 x_0- 4 x_3x_1^3.\]
Equivalently, $(D_A=0)$ is the dual variety of the projectively embedded toric
variety parametrized by monomials with exponents in $A$, i.e., the well known twisted cubic.

Consider the following choice of matrix $B$: 
$$B \, = \, \left(
\begin{array}{rr}
1 & 2 \\
-2 & -3 \\
1 & 0 \\
0 & 1
\end{array}
\right).$$
Note that $g_B=1$, so that the columns of $B$ are a basis of the
integer kernel of $A$.
Calling $l_0(u,v):=u+2 v$\ , $l_1(u,v):=-2 u-3 v$\ , $l_2(u,v):=u$ and $l_3(u,v):=v$,
the parametrization $\psi_B$ equals
\[\left\{
  \begin{array}{lcl}
y_1 & := & \frac{l_0 l_2}{l_1^2}\\
 y_2 & := & \frac{l_0^2 l_3}{l_1^3}\\    
  \end{array}
\right.
\]
Its closed image is the hypersurface $S_B = \{ \Delta_B = 0 \}$, where
\[ \Delta_B (y_1, y_2) = -4 y_2-27 y_2^2+y_1^2+18 y_2 y_1-4 y_1^3,
\]
which can be computed  in \Maple as 
\[-normal(resultant( l_1^2*y_1-l_0*l_2, l_1^3*y_2-l_0^2*l_3, u)\; / \;v^6)\]
(the factor $v^6$ appears because \Maple computes an affine resultant instead of an
homogeneous one with respect to $(u:v)$). 
Equivalently, it could be computed as the dehomogenization
\[\Delta_B(y_1, y_2) = D_A( 1,1,y_1, y_2).\]  
Conversely, up to a monomial $D_A$ equals 
$\Delta_B (x_0 x_2/ x_1^2, x_0^2 x_3 / x_1^3)$.
The associated logarithmic  Gauss map equals
$\gamma_B = (\gamma_1,\gamma_2)$, with
\[\gamma_1:= 12 y_1^3-2 y_1^2-18 y_1 y_2 \quad , \quad
\gamma_2:= -4 y_2-54 y_2^2 +18 y_1 y_2.\] If we compute
$({z_1}{(u,v)},{z_2}{(u,v)}):=(\gamma_B \circ \psi_B)(u,v)$, we obtain:
\begin{eqnarray*}
  z_1 & = & 4 (u+2 v)^2 u (u^3+9 u^2 v+27 u v^2+27 v^3)/(2 u+3 v)^6 \; ,\\
  z_2 & = & 4 (u+2 v)^2 v (u^3+9 u^2 v+27 u v^2+27 v^3)/(2 u+3 v)^6. \\
\end{eqnarray*}
Note that we do not recover $(u,v)$  but it only holds that $z_1/z_2
= u/v$, or $(z_1 : z_2 ) = (u:v)$. 

Consider now the following matrix $C$:
$$C \, = \, \left(
\begin{array}{rr}
1 & 2 \\
0 & -3 \\
-3 & 0 \\
2 & 1
\end{array}
\right).$$
Note that $A \cdot C \, = \, 0$ but that $g_C=3$. In fact, $C \, =  \, B \cdot M$, where
$M$ is the square matrix:
$$ M \, = \, \left(
\begin{array}{rr}
-3 & 0 \\
2 & 1
\end{array}
\right). $$
Calling $L_0(u,v):=u+2v$\ , $L_1(u,v):=-3 v$\ , $L_2(u,v):=- 3 u$ and $L_3(u,v):= 2 u + v$, 
the parametrization $\psi_C$ equals
\[\left\{
  \begin{array}{lcl}
y_1 & := & \frac{L_0 L_3^2}{L_2^3}\\
 y_2 & := & \frac{L_0^2 L_3}{L_1^3} \\    
  \end{array}
\right.
\]
Its closed image is the hypersurface $S_C = \{ \Delta_C = 0 \}$, where
$$
 \Delta_C (y_1, y_2) \,  = \, -1296  y_2  y_1^3-8748  y_2^2  y_1^3-19683  y_2^3  y_1^3+y_2  y_1+4698  y_2^2  y_1^2 $$
$$  -64  y_1^3-64  y_2^3+24  y_2^2  y_1+24  y_2  y_1^2-1296  y_2^3  y_1-8748  y_2^3  y_1^2.
$$
The corresponding logarithmic Gauss map $\gamma_C (y)$ 
composed with $\psi_C (u,v)$ gives back
$ (r(u,v) \cdot u, r(u,v) \cdot v)$, where $r$ denotes the rational function
$ r(u,v) \, = \, -4/729  (u+2  v)^3  (2  u+v)^3  (16  v^4+84  u  v^3+143  u^2  v^2+84  u^3  v+16  u^4)  (v+u)^3/(u^6  v^6)$.
We will study in Section~\ref{sec:three} the precise relation between the irreducible polynomials
$\Delta_B$ and $\Delta_C$.
\end{example}

\begin{proof}[Proof of the ``if'' part in Theorem~\ref{th:Kapranov}]
Let $C$ be a regular non defective $n \times m$ integer matrix, a point $\lambda \in (\CC^*)^m$
in the torus, and consider the map $\psi'_C:=\lambda \psi_C$. We need to show that the
logarithmic Gauss map is its birational inverse. Denote by
$\Delta$ an irreducible equation of its closed image.
The principal observation is that the Jacobian matrix of $\log( \psi'_C)$ is symmetric
since
$$ \frac{\partial} {\partial u_k} \log(  (\psi'_C)_j) \, = \, \sum_{i=1}^n 
\frac{c_{i,k} \ c_{i,j}}{l_i(u)}.$$
Moveover, a straightforward computation shows that 
for any point $u$ in the preimage of the torus, the Jacobian matrices
$ J(\psi'_C)$ and $J(\log( \psi'_C)$ have the same rank  since
$J(\log(\psi'_C) = J(\psi'_C) \cdot D$, where $D$ is the diagonal matrix
with diagonal entries the multiplicative inverses of the coordinates of $\psi'_C$.
This rank is equal to $m-1$ by our hyphotesis that $C$ is non defective.
Now, on one side, implicit partial differentiation of the equality $\Delta( \psi'_C (u)) =0$
implies that the vector $\gamma_C(y)$ lies in the kernel of  the transposed Jacobian
matrix $J (\log ( \psi'_C))^t$ for any $y$ in the image of $\psi'_C$. On the other side,
since the coordinates of $\psi'_C$ are homogeneous forms of degree $0$, it follows
from Euler's formula applied to the coordinates of $\log(\psi'_C))$
that any point $u$ in the preimage of the torus lies in the
kernel of $J(\log(\psi'_C)(u)$. Then, $u$ is proportional to $\gamma_C(\psi'_C(u))$, when
this vector is non zero.
\end{proof}

\section{Monomial changes of coordinates and factorizations} \label{sec:two}

Given any matrix $C \in \ZZ^{n \times m}$ with $g_C \not=0$, we proved that the
parame\-tri\-zation $\psi_C$ is birational if  $C$ is non defective (i.e. when
$S_C$ has codimension $1$). Indeed, the converse is trivially true.
The defectiveness condition can be checked by computing the generic rank of the
Jacobian matrix $J(\psi_C)(u)$. The $(m-1) \times (m-1)$ minors of this matrix
are rational functions, and with probability $1$ at least one of them will
be non zero at a generic point $u$ when $C$ is non defective. The tropical approach
in \cite{dfs} provides another algorithm to check this property.

It is obvious that given a regular non defective matrix $C \in \ZZ^{n \times m}$,
we can replace all row vectors in $C$ lying in the same one-dimensional flat
$\mathcal F$ by their sum, without essentially changing the coordinates of the parametrization
$\psi_C$ except for constants  (if the sum gives the zero vector, we keep the constants but
we don't keep a zero row). So, we can always start simplifying the problem by
making this replacement. Note however that we could change the gcd, so this
operation could lead from a Gale dual matrix with gcd equal to $1$ to a
matrix with arbitrary (non zero) gcd \cite{cc}.

Assume $C$ and $B$ are non-defective  $n \times m$ regular integer matrices of full rank $m$
such that the lattice $\ZZ B$ generated by the columns of $B$ is the saturated of the lattice
$\ZZ C$. Then, there exists a square matrix $M \in \ZZ^{m \times m}$ with determinant
equal to $\pm g_C$ such that $C \, = \, B \cdot M$, as in Example~\ref{cubic}. The lattice
ideal $I(\ZZ C)$ (in $n$ variables) is radical with $|g_C|$ primary components, which correspond to
torus translates of the toric variety defined by the lattice ideal $I(\ZZ B)$ \cite{es}.
We will see in Theorem~\ref{prop:fact} 
how this is reflected in the precise relation between the irreducible $m$-variate
polynomials $\Delta_B$ and $\Delta_C$.

\smallskip

We begin by recalling Lemma~6.3.1 in \cite{curran}, which shows the  relation
between the Horn-Kapranov parametrizations associated to two regular integer $n \times m$ 
matrices $C_1, C_2$ with $C_1 = C_2 \cdot M$ for a square matrix $M$. 

We associate to $M = (M_{ij}) \in \ZZ^{m\times m}$ a linear map $\Lambda_M: \PP^{m-1} \rightarrow
\PP^{m-1}$:
\begin{equation} \label{eq:lambda}
\Lambda_M(u)= (\sum_{j=1}^m{M_{1,j}u_j}:
  \ldots : \sum_{j=1}^m{M_{m,j}u_j})= M \cdot u^t \; ,
\end{equation}
and, denoting by $M^{(1)}, \dots, M^{(m)}$ the column vectors of $M$, the (multiplicative)
monomial map $\alpha_M:{(\CC^*)}^m \rightarrow {(\CC^*)}^m$:
\begin{equation}\label{eq:p}
\alpha_M(y)
=(\prod_{i=1}^m{y_i^{M_{i,1}}}, \ldots,
\prod_{i=1}^m{y_i^{M_{i,m}}})
=(y^{M^{(1)}}, \ldots,
y^{M^{(m)}})\;.
\end{equation}

For any given $m \times m$
matrices $M_1, M_2$, it clearly holds that 
$\Lambda_{ M_1 \cdot M_2} =  \Lambda_{M_1} \circ \Lambda_{M_2}$ and 
\begin{equation}\label{eq:product}
\alpha_{M_1 \cdot M_2} =  \alpha_{M_2} \circ \alpha_{M_1}\; .
\end{equation}

\begin{lemma} \label{lem:square}
 The following diagram 
\[
\xymatrix
{
{\PP^{m-1}} \ar@{-->}[d]_{\psi_{C_1}} \ar@{->}[r]^{\Lambda_M}
 & {\PP^{m-1}} \ar@{-->}[d]^{\psi_{C_2}} \\
 {{(\CC^*)}^m} & {{(\CC^*)}^m} \ar@{->}[l]^{\alpha_M}
}
\]
is commutative.
\end{lemma}

The proof of Lemma~\ref{lem:square} is a straightforward verification.

\begin{remark} \label{rmk:cazabobos}
Assume $C_2$ has gcd equal to $1$, and call $C=C_1, B=C_2$. Then, $|\det(M) | =
g_C$. Suppose that
we didn't know Theorem~\ref{th:Kapranov} but instead we  suspected (or proved)
that $\psi_B$ is birational.
From the equality $\psi_C \, = \alpha_M \circ \psi_B \circ \Lambda_M$, where
$\Lambda_M$ is birational and $\alpha_M$ is a $g_C$ to $1$ mapping, one is
tempted to deduce that $\psi_C$ is also a $g_C$ to $1$ mapping. But indeed, we
have already proved that it is birational. The explanation is given in the next lemma.
\end{remark}

\begin{lemma}
\label{pr:p_proper_dash}
 With definitions and notations as in (\ref{eq:lambda}) and (\ref{eq:p}), the restriction of
 the map $\alpha_M$ to the zero set of $\Delta_{C_2}$ defines by corestriction
 a birational map
\[
 \tilde{\alpha}_M={{\alpha_M}_|}_{(\Delta_{C_2} = 0) }: ( \Delta_{C_2} = 0 )
  \dashrightarrow ( \Delta_{C_1} = 0 ),
\]
where as before $C_1 = C_2 .M$ and $C_1, C_2$ are non defective regular integer matrices.
\end{lemma}

\begin{proof} 
The image of $\psi_{C_2}$ is dense in $(\Delta_{C_2})=0$, and so the
image of $\psi_{C_2} \circ \Lambda_M$ is also dense. For any point $y$
of the form $y = \psi_{C_2} (\Lambda_M(u))$, its image by $\alpha_M$
equals the point $\psi_{C_1}(u)$. Therefore, it lies in $(\Delta_{C_1} =0)$.
Thus, we have a rational map $\tilde{\alpha}_M:  ( \Delta_{C_2} = 0 )
 \dashrightarrow ( \Delta_{C_1} = 0 )$
which has to be 1--1 by Theorem~\ref{th:Kapranov}
  applied to both $\psi_{C_2}$ and $\psi_{C_1}$.
\end{proof}

We now present the relation between  $\Delta_{C_1}$ and $\Delta_{C_2}$. Note that since the
corresponding Horn-Kapranov parametrizations are  given by  rational forms with rational coefficients,
we can assume that these polynomials have integer coefficients and content $1$. They are thus defined
up to sign, but we will usually omit this sign in the notation.

Given an integer square matrix $M$ of size $m$, denote by $G_M$ the multiplicative
group
\begin{equation}\label{eq:G}
G_M:= \, \{ \varepsilon \in (\CC^*)^m \, : \, \alpha_M(\varepsilon) = (1, \dots, 1) \}
\end{equation}
with the induced coordinatewise multiplication. We then have:

\begin{theorem}\label{prop:fact}
Let $C_1, C_2$ are non defective $ n\times m$ regular integer matrices such that 
$C_1 = C_2 .M$. There exists $v$ in the lattice $\ZZ M$ spanned by the
columns of $M$ (or equivalently, verifying $\varepsilon^v =1$ for all 
$\varepsilon \in G_M$)  such that 
\begin{equation} \label{eq:factor}
   \Delta_{C_1} \circ \alpha_M (y) =  y^v  \prod_{\varepsilon \in G_M}{\Delta_{C_2}
  (\varepsilon \cdot y)}\;. 
\end{equation}
\end{theorem}

Before giving the proof, we revisit Example~\ref{cubic}.

\begin{example}[Example~\ref{cubic} cont.]
With the notations of Example~\ref{cubic}, the group $G_M$ consists of those
$(\varepsilon_1, \varepsilon_2) \in (\CC^*)^2$ such that 
$\varepsilon_1^3 =1$, $\varepsilon_2 = 1$ and so
we have
\[
\Delta_C(\frac{y_2^2}{y_1^3}, y_2) \, = \, \frac{y_2^3}{y_1^9} \, \prod_{\varepsilon^3=1}
\ \Delta_B( \varepsilon y_1, y_2).
\]
Of course, we can move the factor $y_1^9$ to the left hand side to 
produce an equality in the polynomial ring $\CC[y_1, y_2]$ but this is really an
equality over the Laurent polynomial ring $\CC[y_1^{\pm 1},y_2^{\pm 1}]$.
The exponent vector $v =(-9,3)$ equals  $3$ times the difference of the columns of $M$.
\end{example}

\smallskip

\begin{proof}
Given any point $y \in (\CC^*)^m$ in the image of $\psi_{C_2}$ (thus, in a dense subset of
$(\Delta_{C_2} =0)$), using that $\Lambda_M$ is an
isomorphism and Lemma~\ref{lem:square}, we can write $ y = \psi_{C_2} (\Lambda_M(u))$ for some
$u \in \PP^{m-1}$  and so $\alpha_M(y) \in (\Delta_{C_1} =0)$. For any $\varepsilon \in G_M$
we have that $\alpha_M(y) = \alpha_M(\varepsilon \cdot y)$ and therefore it also holds that
$\Delta_{C_1}(\varepsilon \cdot y) = 0$. Reciprocally,
pick any point $y \in (\CC^*)^m$ such that $\alpha_M(y)$ lies in the image of $\psi_{C_1}$.
Then, there exists $u \in \PP^{m-1}$ such that $\alpha_M(\psi_{C_2} (\Lambda_M(u)) )=
\psi_{C_1}(u) = \alpha_M(y)$. Therefore, there exists $\varepsilon \in G_M$ such that
$y = \varepsilon^{-1} \psi_{C_2}(\Lambda_M(u))$, and then 
$\Delta_{C_2} (\varepsilon \cdot y) = 0$. By  density and properness arguments, we deduce that
that
\[
(\Delta_{C_1} \circ \alpha_M (y) =0) \cap (\CC^*)^m \, = \,
\bigcup\limits_{\varepsilon \in G_m} (\Delta_{C_2} (\varepsilon \cdot
y ) = 0) \cap (\CC^*)^m
\]

Observe now that, as a consequence of Lemma~\ref{pr:p_proper_dash}, the irreducible polynomials
$\Delta_{C_2}( \varepsilon \cdot y)$ with $\varepsilon$ varying in $G_M$, are pairwise coprime. 
In fact,  if $\Delta_{C_2}( \varepsilon \cdot y)$ is proportional to
$\Delta_{C_2}( \varepsilon' \cdot y)$, with $\varepsilon, \varepsilon' \in G_M$, then
writing $\delta = \varepsilon' \cdot \varepsilon^{-1}$, we have that a point
$y \in (\Delta_{C_2} =0) \cap (\CC^*)^m$ if and only if the point $\delta \cdot y
\in (\Delta_{C_2} =0) \cap (\CC^*)^m$. Since $\tilde{\alpha}_M$ is birational, 
we deduce that $\delta$ is the unit element in $G_M$, i.e. that $\varepsilon =
\varepsilon'$.
Then, we deduce from the Nullstellensatz that there exist positive integers
$n_\varepsilon$, and a unit $ q y^v$ (where $q$ is a constant and $v \in \ZZ^m$)
in the Laurent polynomial ring $\CC[y_1^{\pm 1}, \dots,
y_m^{\pm 1}]$ such that
\begin{equation} \label{eq:nepsilon}
\Delta_{C_1} (\alpha_M(y)) \, = \, q y^v \prod_{\varepsilon \in G_M} \Delta_{C_2} (\varepsilon \cdot y)^{n_\epsilon}.
\end{equation}
Substituting $y \mapsto \delta \cdot y$ for any $\delta \in G_M$ in
the above factorization we get
\[
\Delta_{C_1} (\alpha_M( y)) \, = \, 
\Delta_{C_1} (\alpha_M(\delta \cdot y)) \, = \, q \delta^v  y^v \prod_{\varepsilon \in G_M} \Delta_{C_2} (\varepsilon \cdot y)^{n_{\epsilon \cdot \delta^{-1}}}.
\]
By  uniqueness of the irreducible
factorization, it follows that all $n_\varepsilon$ are  equal to some $N\in \NN$ and that moreover
$\delta^v = (1, \dots, 1)$ for all $\delta \in G_M$. 

It is clear that this last property holds whenever $v$ lies in $\ZZ M$. To prove the converse, assume that
$\delta^v=1$ for all $\delta \in G_M$. To see that $v  \in \ZZ M$ we write
$M$ in its Smith normal form:
\begin{equation}
  \label{eq:smith}
M \,  = \, U \cdot
\underbrace{\left(
\begin{array}{ccc}
d_1 &  &0\\
 & \ddots  & \\
0 &  &d_m
\end{array}
\right)}_{=D}\cdot V,
\end{equation}
where $U,V \in \ZZ^{m\times m}$ are invertible over $\ZZ$ and $d_1\ | \ 
d_2\ |\ \ldots \ | \ d_m$  in $\ZZ$. By the composition formula~(\ref{eq:product}),
it is enough to prove a similar result for $U, D$ and $V$. This is obvious for
$U$ and $V$ since $\ZZ U = \ZZ V = \ZZ^m$ and it can be easily proved for the
diagonal matrix $D$ by choosing elements in $G_D$  of the form $(1, \dots, w_k, \dots, 1)$
with $w_k$ a primitive root of unity of order $d_k$, for all $k=1,\dots, m$.

Assume that all $n_\varepsilon = N > 1$ and differentiate both sides of  equation
~(\ref{eq:nepsilon}). Since the Jacobian of
$\alpha_M$ is invertible at any point of the torus, we deduce from
the Chain rule and the fact that  $\tilde{\alpha}_M$ is birational,
that the Jacobian of $\Delta_{C_1}$  vanishes 
along $(\Delta_{C_1} = 0)$. 
But this contradicts the irreducibility of
$\Delta_{C_1}$, so $N$ must equal 1.

Finally, we show that $q = \pm 1$. As $\Delta_{C_1}$ has content $1$,
we need to show that the product 
\[P(y):= \, \prod_{\varepsilon \in G_M} \Delta_{C_2} (\varepsilon \cdot y)
\]
has integer coefficients and content $1$ too. First, note that if we write $M$ in its Smith Normal
Form (\ref{eq:smith}), it is enough by formula (\ref{eq:product}) to prove the result for
the factors $U,V,D$. This is obvious for $U$ and $V$. We can further decompose $D$ as a
product of diagonal matrices, all of whose diagonal entries are equal to $1$ except for a single entry which is a prime $p$. So, assume that $D= (d_{ij})$ is the diagonal matrix with $d_{ii} =1$
for all $i=1,\dots, m-1,$ and $d_{mm} = p$. But then, the coefficients of $P$ are symmetric
polynomials with coefficients in $\ZZ[s_1, \dots, s_p]$, where $s_1, \dots, s_p$ are the elementary
symmetric functions on the $p$-th roots of unity. Since all $s_i$ equal either $1, 0 $ or $-1$,
we deduce that $P \in \ZZ[y_1, \dots, y_{m-1}][y_m]$.  In fact, there exist a polynomial $Q$ with {\sl the same coefficients\/}, such that $ P (y_1, \dots, y_m) = Q(y_1, \dots, y_{m-1}, y_m^p)$. Moreover,
for every fixed values of the first $m-1$ coordinates, the roots of $Q$ are the $p$-th powers of the roots of $\Delta_{C_2}$ in  the last  variable $y_m$, so that  one can trace recursively the relation between the coefficients of $Q$ and the coefficients of $\Delta_{C_2}$ (which has content $1$) to deduce that the gcd of the coefficients of $Q$ is also equal to $1$. 

A different argument to prove that $q =\pm 1$ is the following. Assume
again that we have a diagonal matrix with all diagonal entries equal
to $1$, except for a single entry which equals a prime number $p$.
Fix a primitive $p$-th root of unity $w$ and take any monomial degree
ordering $\prec$. Call $b_{\gamma}y^\gamma$ the leading term of
$\Delta_{C_2}(y)$. Then, the leading term of $P$ is $b_{\gamma}^g
w^{(\gamma_m\sum_{i=1}^g{i})}y^{g\gamma}= \pm 1 \, b_\gamma^g
y^{g\gamma}$.  Since $\Delta_{C_1} \in \ZZ[y]$, then $q \in \QQ$. But
as all the coefficients of $P$ lie in $\ZZ(w)$, and $q \in \QQ$ we
deduce that $P \in \ZZ[y]$. In fact, $q=1/s$ where $s=${cont}($P$).
Assume $s \neq \pm 1$ and let $a$ be a prime dividing $s$.  Suppose
first that $a \neq p$. Since $a$ divides the content of $P(y) \in
\ZZ[y]$ we have that $P = 0$ in the extension field $\ZZ_a(w)$.
Therefore, one of the factors of $P$ must be zero.  But given that the
content of $\Delta_{C_2}$ is $\pm 1$, this cannot happen.  If $a=p$,
first reduce $\Delta_{C_2}$ mod $p$ and then look at its leading
coefficient, which we call $b_\gamma y^\gamma$. Then, the coefficient
of the monomial $y^{g\gamma}$ in $P$ is not divisible by $p$, a
contradiction.
\end{proof}

The moral of the factorization provided by Theorem~\ref{prop:fact} is that
when trying to compute a generalized discriminant $\Delta_C$ it is possible
to compute a ``simpler'' discriminant $\Delta_B$ where $\ZZ C \subseteq \ZZ B$
and then reconstruct $\Delta_C$ from (\ref{eq:factor}). Even if the monomial
$y^v$ is in principle unknown, its function is to clear denominators without
introducing extra monomial factors. 
The group $G_M$ can be computed via the
Smith Normal Form (\ref{eq:smith}) of $M$. One easier way to recover $\Delta_{C_1}$
from $\Delta_{C_2}$ is the following.
Substitute $y = \alpha_{{\rm Adj} (M)} (z)$ in (\ref{eq:factor}) and 
denote $g : = g_{C_2}$. The key point is  that $\alpha_M
\circ \alpha_{{\rm Adj} (M)} (z) = z^g$. Furthermore, since we know that
the exponent $v$ lies on the lattice spanned by the columns of $M$,
 we deduce that the specialized product on the left hand side must be a Laurent 
polynomial in the variables $z_1^g, \ldots, z_m^g$.  So, we only need to
divide all exponents in the resulting polynomial by $g$ to recover $\Delta_{C_1}(z)$.

In fact, this is also useful when dealing with saturated lattices $\ZZ B$, i.e.
when trying to compute $A$-discriminants. 
Given a regular $n \times m$ non defective matrix $B$, we can first look for
a reduced basis of $\ZZ B$ using the LLL-algorithm \cite{lll},
available at any Computer Algebra System, and then put these reduced generators
as the columns of a matrix $B'$. Write $B  \, = \, B' \cdot M$ with $\det(M)=1$.
Then $G_M$ consists of the single element $(1, \dots, 1)$ and  so (\ref{eq:factor})
reduces to the equality 
\[ \Delta_B (y)   \, = \, y^v  \Delta_{B'} (\alpha_{M^{-1}}(y)) \; ,\]
which can in fact be easily proved since the homogenizations of both
discriminants give the same discriminant $D_A$ (Here, $B$ and $B'$ are Gale duals
of $A$).
Since the coefficients in the parametrization $\psi_{B'}$ are smaller,
the implicit equation $\Delta_{B'}$ can be obtained using standard elimination
techniques in cases in which the systems would crash when trying to compute
$\Delta_B$ directly.

Here is a simple example.

\begin{example} \label{ex:dosbes}[Example~\ref{cubic} cont.]
Recall that for the matrix $B$ we proposed in Example~\ref{cubic}, the inhomogeneous
discriminant had degree $3$. Consider now another choice $B'$ of a Gale dual
of $A$ whose entries are integers with larger absolute value:
\[
B' \, = \, \left(
\begin{array}{rr}
-5 & -3 \\
13 & 8 \\
-11 & -7 \\
3 & 2
\end{array}
\right).
\]
If we try to compute the inhomogeneous discriminant $\Delta_{B'}$ using the resultant
formula, we get a common factor in the coefficients equal to $3^{88}$. After dividing by
this quantity, we recover the following polynomial of degree $16$:
\[
\Delta_{B'}(y) \, = \,-27 y_2^{16}+ 18 y_2^8 y_1^5 - 4y_2^5 y_1^7 - 4 y_2^3 y_1^8 + y_1^{10}.
\]                      
\end{example}

Note that in general, the degree of a dehomogenization of a sparse discriminant $D_A$ 
can be as large as wanted. 

\section{The degree of $\Delta_C$ and the computation of local
multiplicities} \label{sec:three}

Let $C \in \ZZ^{n \times m}$ be, as before, a regular non defective integer matrix with no
zero rows. We set $m=3$ and we assume  w.l.o.g. that the variety of base
points is finite. Base points $p_{\mathcal F} \in {\mathcal Z}$ are indexed by 
basic flats ${\mathcal F}$ as in Section~\ref{sec:one}. 
In this section, we concentrate on the algorithmic computation of the degree of the
generalized homogeneous variety $S_C$,  based on the following well known intersection theory
formula \cite{fulton}
\begin{equation}\label{eq:fulton}
d_C^2 \, = deg(\psi_C) \ deg(S_C) + \sum_{{\mathcal F} \ {\rm basic}} e_{\mathcal F},
\end{equation}
where $e_{\mathcal F}$ denotes the Hilbert-Samuel multiplicity 
of  $p_{\mathcal F}$ \cite{matsumura, BH, slides}.

Since we have an easy formula for $d_C$ in (\ref{eq:d}) and we know that $deg(\psi_C) =1$
by Theorem~\ref{th:Kapranov}, we would need to compute the Hilbert-Samuel multiplicities.
In fact, this is a delicate notion and there is no efficient deterministic
algorithm for the general case (cf. \cite{singularbook,mora} for a Gr\"obner/Standard bases approach). 
We will start by recalling the definition of multiplicity,
together with some known properties.  In particular, there is a probabilistic algorithm,
which reduces the problem to the tractable case of local complete intersection. We will present
however several examples, which show that in general one cannot expect that the 
base points that occur in discriminant parametrizations are local complete intersections,
or even {\em almost local complete intersections} with one more generator. 
It follows that the known algorithms to
find the implicit equation do not work in principle in these cases \cite{Bu1, bcd}. 
One way out would be to compute the Newton polytope of $\Delta_C$ using the results
in \cite{dfs}, and then compute its coefficients via an efficient interpolation.

We show in Proposition~\ref{prop:uniform} 
that when  $C$ is {\em uniform}, i.e. when all its maximal minors are non zero, 
one can use the combinatorial algorithm from \cite{vol_alg} to compute the Hilbert-Samuel
multiplicities.  We then turn things upside-down in Corollary~\ref{cor:sparse}
to compute the dimension of the local
vector space at the origin of $d$ sparse polynomials in $d$ variables with generic
coefficients.

There are several ways to define the algebraic multiplicity of a base
point.  Our definition follows \cite{matsumura}. We will state the results
for the case of dimension two, but they hold with the obvious changes
in any dimension. We refer to~\cite{tesis} for more details on the definitions 
and examples in this section. 

\smallskip

Let $p = p_{\mathcal F} \in {\mathcal Z}$. Consider the Noetherian local ring
$A_p:={\mathcal O}_{{\PP^2}, p}$ and the localized base point locus ideal 
$I_p:=\langle f_0, f_1, f_2, f_3 \rangle A_p$.
The Samuel function of $A_p$ with respecto to $I_p$ is defined as:
\[
\chi_{A_p}^{I_p}(r)= l(A_p/I_p^{r+1}) \quad \textrm{for all } r \in  \NN \; ,
\]
were $l(\_)$ is the length function of $I_p$ as an $A_p$-module, that is, the
length of a composition series of the module. Since we are working over the algebraically
closed field $\CC$, this length coincides with the vector space dimension
$\dim_\CC(A_p/I_p^{r+1})$.
The Samuel function is polynomial for large
values of $r$, that is,
there is a polynomial $PS_{A_p}^{\; I_p}(X)$ in $\QQ[X]$ (which takes integer values over $\ZZ$) 
such that we have
$PS_{A_p}^{\; I_p}(r)=\chi_{A_p}^{I_p}(r)$ for $r >>0$. Moreover, this polynomial has
degree $2$ and  its leading coefficient is
$e/2!$ with $e \in \NN_0$. Then, the local multiplicity $e_{\mathcal F}$ of the base point locus
at $p_{\mathcal F}$ is defined as $e_{\mathcal F}:= e$, i.e. $2!$ times the leading
coefficient of the polynomial $PS_{A_p}^{\; I_p}$.

Therefore, 
\[
  e_{\mathcal F} \, = \, \lim_{r \rightarrow \infty}{\frac{\dim_\CC(A_p/I_p^{r+1}) \cdot 2!}{r^2}}\;.
  \]
When the base point $p = p_{\mathcal F}$ is a local complete intersection, i.e., when
the ideal $I_p$ admits two generators, then $e_{\mathcal F}$ is just the vector space dimension
$A_P/ I_p$ of the local quotient.
This dimension can thus be computed algorithmically via a standard basis computation
using a local order $\prec$, and counting the number of monomials not in
$in_\prec(I_p)$. Even if this is not the case,  we always have the following probabilistic
approach to compute the local multiplicity.

Consider the ideal $J_p$ generated by
$2$ generic linear combinations of the $4$ generators:
\[
J_p:=<v_0^{0} f_0+v_1^{0} f_1+ v_2^{0} f_2+v_3^{0} f_3\; ,
\;v_0^{1} f_0+v_1^{1} f_1+ v_2^{1} f_2+v_3^{1} f_3>\; ,
\]
with $v_i^{j} \in \CC$. Then $J_{p}$ is generically a complete
intersection inside $I_p$ (and a {\em reduction ideal} of $I_p$). Thus, with probability $1$,
we can compute
$e_{\mathcal F} =\dim_\CC{ (A_p/J_{p})}$.

As a corollary, we always have the inequality $e_{\mathcal F}  \geq \dim_\CC{ (A_p/J_{p})}$, so
that in any case
$$ deg(\Delta_C) \, \leq \, d_C^2 - \sum_{{\mathcal F} \, {\rm basic}} \dim_\CC (A_{p_{\mathcal F}}/
I_{p_{\mathcal F}}).$$

On the other side when $I$ is a monomial
ideal, there exists a combinatorial way of computing this
multiplicity, as stated in \cite{vol_alg}. We reproduce this result below. 

If $p=(1:0:0)$ is a base point (which we can assume after a
translation) and the localized ideal $I_p$ is monomial,
 we have the following algorithm to compute the
Hilbert-Samuel multiplicity $e_p$ at $p$:

\begin{algorithm}\label{alg:Vol}
 Computation of Hilbert-Samuel Multiplicities for the mono\-mial
  case and $m=3$. 
  \begin{itemize}
  \item Set $x_0 = 1$ and let $\tilde{I}_p$ be the specialization
  of the ideal $I_p$.
  \item Compute the convex hull $\mathcal{C}$ of the exponents of the
  bivariate monomials in $\tilde{I}_p$.
  \item Then: $e_{p} = 2! \cdot Vol(\NN_0^2 \smallsetminus
    \mathcal{C})$ equals the normalized volume of the complement $\mathcal K$ of
    $\mathcal C$ in the first orthant.
  \end{itemize}
\end{algorithm}

Note that in the very simple case in which the ideal is both monomial
and a complete intersection, generated by $\{x_1^{m_1}, x_2^{m_2}\}$, 
the local multiplicity $m_1 \times m_2$ equals both the normalized
volume of the triangle $\mathcal K$ with vertices $(0,0), (m_1,0), (0, m_2)$, which is
the complement in $\NN_0^2$ of the convex hull of the staircase of the
ideal, and the number of standard monomials $\{x_1^{k_1} x_2^{k_2}\, / \,
0 \leq k_i < m_i, i=1,2\}$, which is the dimension of the quotient
by the ideal.

\begin{remark} \label{rmk:novale}
The Volume Algorithm \ref{alg:Vol} does not work for general
ideals. If fact, it might seem reasonable to expect the same algebraic multiplicity for
$I_p$ and for any initial (monomial) ideal $in_\prec(I_p)$ with respect to a local order, 
but this is not in general the case. 
We illustrate this issue in Example~\ref{ex:A}.
\end{remark}

\medskip
\begin{proposition} \label{prop:uniform}
Assume $C \in \ZZ^{n \times 3}$ is uniform.
Then $C$ is non defective and each local ideal $I_p$ becomes a monomial
ideal modulo a linear change of coordinates. So, the degree of the 
generalized discriminant surface $S_C$ can be combinatorially computed
using formula (\ref{eq:fulton}) and Algorithm~\ref{alg:Vol}. 
\end{proposition}

\begin{proof}
The fact that $C$ is non defective follows from \cite{dfs}, and more precisely
from \cite[Section 5]{cc}. Since any flag is generated by only
two linear forms $l_i, l_j$, after dehomogenizing and localizing,
the four generators of each local ring $I_p$ are products of powers of $l_i, l_j$,
i.e. they are monomials in two independent variables $l_i, l_j$. 
\end{proof}

\begin{example}\label{ex:1}
We first give a very simple example to illustrate Proposition~\ref{prop:uniform}.
Let $C$ be the uniform matrix:
 \[
C = \left(
  \begin{array}{rrr}
    2 & 1 & 3\\
-2 & -1 & -2 \\
1 & 1 & 0 \\
-1 & -1 & -1 \\
  \end{array}
\right)\;.
\]
We read from the parametrization $\psi_C$ that
\[
f_1 = l_1^2 l_3 \quad , \quad f_2 = l_1 l_2 l_3 \quad , \quad
f_3 = l_1^3 \quad , \quad f_0 = l_2^2 l_4 \;.
\]
There are two base points: $p = p_{\{1,4\}} = (-2:1:1)$ and $p' = p_{\{1,2\}}
= (1:-2:0)$. The localized ideals are the following monomial ideals
 $I_{p} =\langle l_1^2 , l_1, l_1^3,
l_4 \rangle = \langle l_1, l_4 \rangle$ and $I_{p'} = \langle l_1^2, l_1 l_2, l_1^3, l_2^2\rangle =
\langle l_1^2, l_1 l_2, l_2^2 \rangle$ (in variables $l_1, l_4$ and $l_1, l_2$ and not in
the $u$ variables). The first ideal is moreover a complete intersection.
It is straightforward in this case to compute the multiplicities: $e_p = 1$ and $e_{p'} = 4$, while
the dimension of the local ring at $p'$ equals $3$. Thus, $deg(S_C) \, = \, 3^2 - 1 -4 \, = 4$.
Indeed, in this case $g_C=1$, so that $C$ is a Gale dual of the matrix $A \in \ZZ^{1 \times 4}$
with all four entries equal to $1$. So, the homogeneous $A$-discriminant $D_A(x) = x_1 + x_2 + x_3 + x_4$.
Then, we can obtain $\Delta_B$ by dehomogenizing $D_A$. From the equations $y_1 = \frac{x_1^2 x_3}{x_2^2 x_4},
\, y_2 = \frac{x_1 x_3}{x_2 x_4}, \, y_3 = \frac{x_1^3}{x_2^2 x_4}$ read from the columns of $C$, we get that
\[ D_A(x) \, = \, x_1 \left(1 + \frac{x_2}{x_1} + \frac{x_3}{x_1} + \frac {x_4}{x_1} \right) \, = \,
x_1  \left( 1 + \frac{y_2}{y_1} + \frac{y_1}{y_3} + \frac {y_2^2 y_3}{y_1^2}\right).
\]
Clearing denominators, we get the equation
\[ \Delta_C( y) \, = \, y_1^2 y_3 + y_1 y_2 + y_1^3 + y_2^2 y_3^2,\]
of degree $4$, as predicted. In fact, this same procedure holds for any matrix $B$ of size $(m+1) \times m$
and $g_B=1$. The associated discriminant $\Delta_B$ has $(m+1)$ monomials and all coefficients are
equal to $1$.
\end{example}

We address now a more complicated example, which will help us illustrate
several features.

\begin{example} \label{ex:A}
Consider the matrix
\[
C=
\left(
\begin{array}{rrr}
1 & -1 & 0 \\
1 & -1 & 1 \\
1 & -1 & 0\\
-1 & 2 & 0 \\
-1 & 1 & -2 \\
-1 & 0 & 1\\
\end{array}
\right) \; .
\]
Observe that the first and the third rows of $C$ are identical. We have:
$f_0 =l_1^4 l_2^2 l_5$ , $f_1 =l_4^3 l_5^3 l_6$ ,
$f_2 =l_1^2 l_2^2 l_4 l_6^2$ , $f_3:=l_1^2 l_2 l_4 l_5^2 l_6$.
There are seven base points:
$p_{\{1,3,4\}} =(0:0:1)$~, $p_{\{1,2,3,5\}}= (1:1:0)$~, $p_{\{1,6\}}=(1:1:1)$~,
$p_{\{2,4\}} =(-2:-1:1)$~, $p_{\{2,6\}}:=(1:2:1)$~,
$p_{\{4,5\}} =(-4:-2:1)$, and $p_{\{5,6\}} =(1:3:1)$.

Let's focus on $p = p_{\{1,3,4\}}$. The local ideal equals $I_p \, = \,
 \langle l_1^4\,;\,l_4^3\,;\,l_1^2 l_4 \rangle$, or changing the name of
 the variables, $I_p \, = \, \langle x_1^4, x_2^3, x_1^2 x_2  \rangle$.
By the volume formula in Algorithm~\ref{alg:Vol} we get $e_p = 10$ (see
Figure~\ref{fig:m10}). The dimension of the local quotient by $I_p$ equals $ 8 < 10$.

On the other side, if we write the linear forms $l_i$ in the
affine coordinates $(u_1, u_2)$, we look at the generators of $I_p$
in the polynomial ring $\CC[u_1, u_2]$. If we consider the local order $\prec =
ds$ (with $u_2 \prec u_1$) in \Singular \cite{singular}, we get the following initial ideal:
\[in_\prec(I_p) \, = \, \langle u_1^3, u_1^2 u_2, u_1 u_2^3, u_2^4 \rangle,
\]
which, by the same volume algorithm as above, has multiplicity $11 > 10$. See
Figure~\ref{fig:m11}.

\begin{figure}
\centering
\includegraphics{dibujitosTesis.1}
    \caption{Region corresponding to the ideal $I_p$}
\label{fig:m10}
\centering
  \includegraphics{dibujitosTesis.2}
   \caption{Region corresponding to the ideal $in_\prec(I_p)$}
\label{fig:m11}
\end{figure}

The implicit equation $\Delta_C(y)$ can in this case be easily computed with
\Singular by Gr\"obner bases methods and we get (up to sign) the following polynomial
of degree $13$:
\[
\begin{array}{l}
-8 y_1^4 y_2^4 y_3^2+3 y_1^2 y_2^2 y_3+16 y_3^3+1000 y_1^3 y_2^2 y_3^2+3 y_1
y_2 y_3
+y_1^3 \\
 y_2^3 y_3+y_3+3125 y_1^4 y_2^2 y_3^2+27 y_1^2 y_2+16 y_1^5 y_2^5 y_3^3-
225 y_1^2 y_2 y_3-\\
-225 y_1^3 y_2^2 y_3+500 y_1^2 y_3^2 y_2+500 y_1^4 y_2^3 y_3^2
+160 y_1^2 y_2^2 y_3^3
+80 \\
 y_1^4 y_2^4 y_3^3-48 y_1^2 y_2^2 y_3^2-32 y_1^3 y_2^3 y_3^2
+160 y_1^3 y_2^3 y_3^3-32 y_1 y_2
 y_3^2-\\
8 y_3^2+80 y_3^3 y_1 y_2 \;.
\end{array}
\]
\end{example}

\begin{example} \label{ex:3}
One could try to reduce the computation of multiplicities to the case of monomial
ideals in the following way. Consider for instance the matrix
\[
C:=
\left(
\begin{array}{rrr}
  1 & 1 & 2\\
-1 & 0 & 1\\
0 & 1 & 3 \\
0 & -1 & -2\\
0 & -1 & -4\\
\end{array}
\right) \;,
\]
and the base point $p = p_{\{2,3,4\}}$. Calling $x_1 = l_4 (1, u_2, u_2), \, x_2 = l_5(1, u_1, u_2)$ 
we have
that $2 l_3 = - (l_4 + l_5)$, so the local ideal $I_p$ is generated by
$$I_p \, =  \langle (x_1+x_2) x_1 x_2^3, (x_1+x_2)^3, x_1^2 x_2^4 \rangle.$$ 
This is not a monomial ideal in these coordinates, but after the linear change
$x_1 = u + 3 v, \, x_2 = - u + v$, it becomes a monomial ideal with the
following generators obtained after easy algebraic manipulations: $\langle u^4 v, u^6, v^3 \rangle$,
and we can compute its multiplicity $e_p = 18$ by means of Algorithm~\ref{alg:Vol}.
However, this case is very special and there doesn't seem to be a general pattern
about when such a change of coordinates is possible. 

For instance, it is possible to prove that
the local ideal at the base point $p_{\{1,4,5\}}$ of the parametrization $\psi_C$
associated to the matrix
\[ C \, = \,
\left(
\begin{array}{rrr}
  1 & 1 & 3\\
1 & 0 & 2\\
0 & 1 & 1 \\
-2 & -2 & 0\\
0 & 0 & -6\\
\end{array}
\right),
\]
cannot be transformed into a monomial ideal by a linear change of coordinates.
\end{example}

\begin{example}
Consider the matrix $C$:
 \[
 C = \left(
   \begin{array}{rrr}
     1 & -7 & -6\\
 -1 & 4 & 3 \\
 1 & 0 & 4 \\
 0 & 1 & -1 \\
 -1 & 2 & 0
   \end{array}
 \right)\;
 \]
and the base point $p = p _{\{1,2\}} =  (-1: -1: 1)$. Calling
$ x = l_1,  \  y = l_2$, the local ideal at $p$ equals
$ I_p = \langle x^8 , y^ 5 , xy^2 , x^7y \rangle$,
which is not an almost complete intersection.
\end{example}

We end with an application to unmixed sparse polynomial systems.
Fix an exponent set $A = \{ \alpha_1, \dots, \alpha_r\} \subseteq \NN_0^d$,
with $r\geq d$, and consider $d$ generic polynomials $F_1, \dots, F_d$
with exponents in $A$ and coefficients in $\CC$:
\[
F_i (x) = \sum_{j=1}^r{c^i_j \, x^{\alpha_j}}\;, \qquad i =1, \ldots, d \quad ;
\quad  x = (x_1, \ldots, x_d)
\]
with ${(c^1_j)}_j, \ldots, {(c^d_j)}_j\in \CC^r$ generic. Then, by Bernstein's
theorem, the total number of  common roots in the torus $(\CC^*)^d$
equals the normalized volume of $A$. Assume moreover than $\alpha_i = \lambda_i \, e_i$
for all $i=1, \dots, d$, i.e. that a monomial which is a pure positive power of 
each of the variable occurs in the  polynomials $F_1, \dots, F_d$. 
Using the previous results, we are
able to compute geometrically their multiplicity at the origin (see also \cite[Chapter 5, 
\S~2.E]{gkzbook} for a general version of this result).

\begin{corollary}\label{cor:sparse}
  Let $A= \{\alpha_1, \ldots, \alpha_r\} \subseteq \NN_0^d$ such that
 $\alpha_i =\lambda_i e_i$,  $i=1, \ldots,
  d$, where $\lambda_i \in \NN$. Given generic sparse polynomials $F_1, \dots,
  F_d$ with exponents in $A$, their multiplicity at the origin
  $$ e_{0} \, = \dim_\CC {\left(\frac{\CC[x_1, \ldots, x_d]}{\langle F_1, \ldots,
    F_d \rangle}\right)}_{0}$$ 
coincides with the normalized volume of the complement ${\mathcal K}_A$ in the
first orthant of the convex hull of $\{\alpha_1 + \RR_{\geq 0 }^{\ d}\} \cup \ldots \cup \{\alpha_r +
\RR_{\geq 0 }^{\ d}\}$.
\end{corollary}

\begin{proof}
By our hypothesis about $A$, it follows that the monomial ideal $I$ generated by
$\{x^{\alpha_1} , \dots, x^{\alpha_r}\}$ is supported at the origin $0 \in \CC^d$.
We identify it with its localization $I_0 = I  \ \CC[x_1, \dots, x_d]_0$.
Note that the localized ideal $J_0:= \langle F_1, \dots, F_d\rangle_0$ is a generic
 complete intersection inside the zero dimensional ideal $I_0$.
By the probabilistic algorithm (in dimension $d$) 
for the computation of the local multiplicity $e$ of
$I_0$, we know that $e$ equals the local multiplicity of the reduction ideal $J_0$
of $I_0$, and then $e = e_0$. Finally, $e$ coincides with the normalized volume
of ${\mathcal K}_A$ by  the corresponding version of Algorithm~\ref{alg:Vol} in dimension
$d$. 
\end{proof}

\subsection*{Acknowledgments}

We thank the Institute for Mathematics and its Applications (IMA) in Minneapolis, 
USA, where this work was completed. We are also grateful to Eduardo Cattani,
 Teresa Krick and Rafael Villarreal for very useful discussions.

\end{document}